\documentclass[english,11pt]{article}
\usepackage[latin1]{inputenc}
\usepackage{amsmath,amsthm,amssymb,babel}

\newtheorem{teo}{Theorem}
\newtheorem{defin}{Definition}

\def\eq#1{(\ref{#1})}
\def\neweq#1{\begin{equation}\label{#1}}
\def\endeq{\end{equation}}
\def\weak{\rightharpoonup}
\def\eps{\varepsilon}
\def\ep{\varepsilon}

\def\phi{\varphi}
\def\supp{{\rm supp}}
\def\C{{C_c^\infty}}
\def\di{\displaystyle}
\def\ri{\rightarrow}
\def\RR{{\mathbb R} }

\def\QQ{{\mathbb Q} }

\def\r2{{\mathbb R}^{2}}

\begin{document}
\title{\sc Entire solutions of multivalued nonlinear Schr\"odinger equations in Sobolev spaces 
with variable exponent}
\author{\sc Teodora-Liliana Dinu \\ \small Department of Mathematics, ``Fra\c tii Buze\c sti" 
College, Bd. \c Stirbei--Vod\u a No. 5, 200352 Craiova, Romania\\ \small 
E-mail: {\tt tldinu@gmail.com}}
  \date{}
      \maketitle
	  
\bigskip
{\small {\bf Abstract.} We establish the existence of an entire solution for a class of stationary 
Schr\"odinger equations with subcritical discontinuous nonlinearity and lower bounded potential 
that blows-up at infinity. The abstract framework is related to Lebesgue--Sobolev spaces with 
variable exponent. The proof is based on the critical point theory in the sense of Clarke and we 
apply Chang's version of the Mountain Pass Lemma without the Palais--Smale condition for locally 
Lipschitz functionals. Our result generalizes in a nonsmooth framework a result of Rabinowitz 
\cite{rabi} on the existence of  ground-state solutions of the nonlinear Schr\"odinger equation.


}
\normalsize

\section{Introduction and auxiliary results}
The Schr\"odinger equation plays the role of Newton's laws and conservation of energy in classical 
mechanics, that is, it predicts the future behaviour of a dynamic system. The linear form of 
Schr\"odinger's equation is
$$\Delta\psi +\frac{8\pi^2m}{\hbar^2}\,\left(E(x)-V(x)\right)\psi =0\,,$$
where $\psi$ is the Schr\"odinger wave function, $m$ is the mass, $\hbar$ denotes Planck's 
constant, $E$ is the energy, and $V$ stands for the potential energy. The structure of the 
nonlinear Schr\"odinger equation is much more complicated. This  equation  
is a prototypical 
dispersive nonlinear
partial differential equation that has been central for
almost four decades now to a variety of areas in Mathematical Physics.
The relevant fields of
application may vary from optics and propagation of the  electric field in optical fibers
(Hasegawa and Kodama \cite{hase}, Malomed \cite{malo}),
to the self-focusing and 
collapse of Langmuir waves in plasma physics (Zakharov \cite{zakh})
and the behaviour of deep water waves and 
freak waves (the so-called rogue waves)
in the ocean (Benjamin and Feir \cite{benj} and Onorato, 
Osborne, Serio and Bertone \cite{onor}). 
The nonlinear Schr\"odinger equation also describes various phenomena arising in: 
self-channelling of a high-power
ultra-short laser in matter, in the theory of Heisenberg ferromagnets and magnons,
in dissipative quantum mechanics, in condensed matter theory, in plasma physics
(e.g., the Kurihara superfluid film equation). We refer to Ablowitz, Prinari and 
Trubatch \cite{abl}, Grosse and Martin \cite{grosse}
and Sulem \cite{sulem} for a modern overview, including applications.

Consider the model problem
\begin{equation}\label{oheq}
i\hbar\psi_t=-\frac{\hbar^2}{2m}\,\Delta\psi+V(x)\psi-\gamma |\psi|^{p-1}\psi\qquad
\mbox{in $\RR^N$ ($N\geq 2$)}\,,\end{equation}
where $p<2N/(N-2)$ if $N\geq 3$ and $p<+\infty$ if $N=2$. 
In physical problems,
a cubic nonlinearity corresponding to $p = 3$ is common; in this case \eq{oheq} is called the
Gross-Pitaevskii equation.
In the study of Eq.~\eq{oheq}, Oh \cite{oh} supposed that the potential $V$ is bounded
and possesses a non-degenerate critical point at $x=0$. More precisely, it is assumed that
$V$ belongs to the class ($V_a$) (for some $a$) introduced in Kato \cite{kato}. 
Taking $\gamma>0$ and $\hbar >0$
sufficiently small and using a Lyapunov-Schmidt type reduction, Oh \cite{oh} proved the existence 
of
a standing wave solution of Problem \eq{oheq}, that is, a solution of the form 
\begin{equation}\label{oheq1}
\psi (x,t)=e^{-iEt/\hbar}u(x)\,.\end{equation}
Note that substituting the ansatz \eq{oheq1} into \eq{oheq} leads to
$$-\frac{\hbar^2}{2}\,\Delta u+\left(V(x)-E\right)u=|u|^{p-1}u\,.$$
The  change of variable $y=\hbar^{-1}x$ (and replacing $y$ by $x$) yields
\begin{equation}\label{oheqoheq}
-\Delta u+2\left(V_\hbar (x)-E\right)u=|u|^{p-1}u\qquad\mbox{in $\RR^N$}\,,\end{equation}
where $V_\hbar (x)=V(\hbar x)$.

If for some $\xi\in\RR^N\setminus\{0\}$, $V(x + s\xi) = V (x)$ for all $s \in \RR$, equation 
\eq{oheq} is invariant
under the Galilean transformation
$$\psi(x,t)\longmapsto \psi(x-\xi t,t)\exp\left(i\xi\cdot x/\hbar -\frac 
12i|\xi|^2t/\hbar\right)\psi(x-\xi t, t)\,.$$
Thus, in this case, standing waves reproduce solitary waves travelling in the direction
of $\xi$. 
In a celebrated paper, Rabinowitz \cite{rabi} 
proved that Equation \eq{oheqoheq} has a ground-state solution (mountain-pass solution) for $\hbar 
>0$ small, under the assumption that $\inf_{x\in\RR^N}V(x)>E$.
After making
a standing wave ansatz, Rabinowitz reduces the problem to that of studying the
semilinear elliptic equation
\begin{equation}\label{phreqq}
-\Delta u +a(x)u=f(x,u)\qquad\mbox{in }\RR^N,
\end{equation}
under suitable conditions on $a$ and assuming that $f$ is smooth, superlinear
and subcritical.

Our purpose in this paper is to study the multivalued version of Equation \eq{phreqq}, but for a 
more general class of differential operators, the so-called $p(x)$--Laplace operators. This 
degenerate quasilinear operator is defined by $\Delta_{p(x)}u:={\rm div}(|\nabla u|^{p(x)-2}\nabla 
u)$ (where $p(x)$ is a certain function whose properties will be stated in what follows) and it 
generalizes the celebrated $p$--Laplace operator
$\Delta_p u:=\mbox{div} (|\nabla u|^{p-2}\nabla u)$, where $p>1$ is a constant.
The $p(x)$--Laplace operator possesses more complicated nonlinearity than the $p$--Laplacian, for 
example, it is inhomogeneous.
We only recall that $\Delta_p$ describes a  variety of phenomena in the nature. For instance, the 
equation
governing the motion of a fluid involves the $p$--Laplace operator. More
exactly, the shear
stress $\vec\tau$ and the velocity gradient $\nabla u$ of the fluid
are related in the manner that
$\vec\tau (x)=r(x)|\nabla u|^{p-2}\nabla u$,
where $p=2$ (resp., $p<2$ or $p>2$) if the fluid is Newtonian (resp.,
pseudoplastic or dilatant). Other applications of the $p$--Laplacian also appear
 in the study of flow through porous media ($p=3/2$),
Nonlinear Elasticity ($p\geq 2$), or Glaciology
 ($1<p\leq 4/3$).  
 
The analysis we develop in this paper is carried out in terms of Clarke's critical point theory 
for locally Lipschitz functionals and in generalized Sobolev spaces. 
That is why we recall in this section some basic facts related to Clarke's generalized gradient 
(see Clarke \cite{3,4} for more details) and Lebesgue-Sobolev  spaces with variable exponent.

 Let $E$ be a real Banach 
space and assume that $I:E\ri\RR$ is a locally Lipschitz functional. Then the Clarke
generalized gradient is defined by
$$\partial I(u)=\{\xi\in E^{*};\, I^{0}(u,v)\geq\langle\xi,v\rangle\,,
\ \mbox{for all}\ v\in E\}\,,$$
where $I^{0}(u,v)$ stands for the directional derivative of $I$ at $u$ in the direction $v$, that 
is,
$$ I^{0}(u,v)=\limsup\limits_{w\to u\atop \lambda\searrow 0}\di\frac{I(w+
\lambda v)-I(w)}{\lambda}\,.$$

Variable
exponent Lebesgue spaces appeared in the literature for the 
first time already in a 1931 article by W.~Orlicz \cite{orl}. In
the years 1950 this study was carried on by Nakano \cite{nak} who
made the first systematic study of spaces with variable exponent. 
Later, the Polish mathematicians investigated the modular function 
spaces (see, e.g., the basic monograph Musielak \cite{M}). Variable exponent Lebesgue 
spaces on the real line have been independently developed by Russian
researchers. In that context we refer to the works of Sharapudinov \cite{sha},
Tsenov \cite{tse} and Zhikov \cite{Z1, Z2}. For deep results
in weighted Sobolev spaces with applications to partial differential equations and nonlinear 
analysis we refer to the excellent monographs by  Drabek, Kufner and Nicolosi \cite{dkn}, 
Hyers,  Isac and Rassias
\cite{isac}, Kufner and Persson \cite{kp}, and Precup \cite{preklu}. 
We also refer to the recent works by
Diening \cite{dien}, Ruzicka \cite{ruz} and Chen, Levine and Rao \cite{ychen}
 for applications of Sobolev spaces with
variable exponent in the study of electrorheological fluids or in
image restoration.

We recall in what follows some definitions and basic properties
of the generalized Lebesgue--Sobolev spaces $L^{p(x)}(\Omega)$
and $W_0^{1,p(x)}(\Omega)$, where $\Omega\subset\RR^N$ is an arbitrary
domain with smooth boundary.

Set
$$C_+(\overline\Omega)=\{h;\;h\in C(\overline\Omega),\;h(x)\geq 2\;{\rm 
for}\;
{\rm all}\;x\in\overline\Omega\}.$$
For any $h\in C_+(\overline\Omega)$ we define
$$h^+=\sup_{x\in\Omega}h(x)\qquad\mbox{and}\qquad h^-=
\inf_{x\in\Omega}h(x).$$
For any $p(x)\in C_+(\overline\Omega)$, we define the variable exponent
Lebesgue space
$$L^{p(x)}(\Omega)=\{u;\ u\ \mbox{is a
 measurable real-valued function such that }
\int_\Omega|u(x)|^{p(x)}\;dx<\infty\}.$$
On this space we define the {\it Luxemburg norm} by 
the formula
$$|u|_{p(x)}=\inf\left\{\mu>0;\;\int_\Omega\left|
\frac{u(x)}{\mu}\right|^{p(x)}\;dx\leq 1\right\}.$$
Variable exponent Lebesgue spaces resemble classical Lebesgue spaces
in many respects: they are Banach spaces \cite[Theorem 2.5]{KR}, the
H\"older inequality holds
\cite[Theorem 2.1]{KR}, they are reflexive if and only if $1 < p^-\leq
p^+<\infty$
\cite[Corollary 2.7]{KR} and continuous functions are dense if $p^+
<\infty$ \cite[Theorem 2.11]{KR}. The inclusion between
Lebesgue spaces also generalizes naturally \cite[Theorem 2.8]{KR}: if 
$0 <
|\Omega|<\infty$
 and $r_1$, $r_2$
are variable exponents so that $r_1(x) \leq r_2(x)$ almost everywhere 
in $\Omega$ then there exists the continuous embedding
$L^{r_2(x)}(\Omega)\hookrightarrow L^{r_1(x)}(\Omega)$, whose norm 
does not exceed $|\Omega|+1$.

We denote by $L^{p^{'}(x)}(\Omega)$ the conjugate space
of $L^{p(x)}(\Omega)$, where $1/p(x)+1/p^{'}(x)=1$. For any
$u\in L^{p(x)}(\Omega)$ and $v\in L^{p^{'}(x)}(\Omega)$ the H\"older
type inequality
\begin{equation}\label{Hol}
\left|\int_\Omega uv\;dx\right|\leq\left(\frac{1}{p^-}+
\frac{1}{{p^{'}}^-}\right)|u|_{p(x)}|v|_{p^{'}(x)}
\end{equation}
holds true.

An important role in manipulating the generalized Lebesgue-Sobolev 
spaces
is played by the {\it modular} of the $L^{p(x)}(\Omega)$ space, which 
is
the mapping
 $\rho_{p(x)}:L^{p(x)}(\Omega)\rightarrow\RR$ defined by
$$\rho_{p(x)}(u)=\int_\Omega|u|^{p(x)}\;dx.$$
If $(u_n)$, $u\in L^{p(x)}(\Omega)$ and $p^+<\infty$ then the 
following relations hold true
\begin{equation}\label{L4}
|u|_{p(x)}>1\;\;\;\Rightarrow\;\;\;|u|_{p(x)}^{p^-}\leq\rho_{p(x)}(u)
\leq|u|_{p(x)}^{p^+}
\end{equation}
\begin{equation}\label{L5}
|u|_{p(x)}<1\;\;\;\Rightarrow\;\;\;|u|_{p(x)}^{p^+}\leq
\rho_{p(x)}(u)\leq|u|_{p(x)}^{p^-}
\end{equation}
\begin{equation}\label{L6}
|u_n-u|_{p(x)}\rightarrow 0\;\;\;\Leftrightarrow\;\;\;\rho_{p(x)}
(u_n-u)\rightarrow 0.
\end{equation}
Spaces with $p^+ =\infty$ have been studied by Edmunds, Lang and 
Nekvinda
\cite{edm}.

Next, we define $W_0^{1,p(x)}(\Omega)$ as the closure of
$C_0^\infty(\Omega)$ under the norm
$$\|u\|_{p(x)}=|\nabla u|_{p(x)}.$$
The space $(W_0^{1,p(x)}(\Omega),\|\cdot\|_{p(x)})$ is a separable 
and reflexive Banach space. We note that if $q\in C_+(\overline
\Omega)$ and $q(x)<p^\star(x)$ for all $x\in\overline\Omega$ then 
the embedding
$W_0^{1,p(x)}(\Omega)\hookrightarrow L^{q(x)}(\Omega)$
is compact (if $\Omega$ is bounded) and continuous (for arbitrary $\Omega$), where 
$p^\star(x)=\frac{Np(x)}{N-p(x)}$
if $p(x)<N$ or $p^\star(x)=+\infty$ if $p(x)\geq N$. We refer to
Edmunds and R\'akosn\'{\i}k
\cite{edm2,edm3}, Fan, Shen and Zhao \cite{FSZ}, Fan and Zhao \cite{FZ1},
and Kov\'a\v cik and R\'akosn\'{\i}k \cite{KR} for further properties of variable exponent
Lebesgue-Sobolev spaces.

\section{The main result}
For any function $h(x,\cdot)\in L^\infty_{{\rm loc}}(\RR)$ we denote by $\underline h$
(resp., $\overline h$) the lower (resp., upper) limit of $h$ in its second variable, that is,
$$\underline h(x,t)=\lim_{\ep\searrow 0}\ {\rm essinf}\ \{h(x,s);\
|t-s| < \ep\}\,;\qquad\overline h(x,t)=\lim_{\ep\searrow 0}\ {\rm esssup}\
\{h(x,s);\ |t-s| < \ep\}\,.$$

Let $a\in L^\infty_{{\rm loc}}(\RR^N)$ be a variable potential such that, for some $a_0>0$,
\neweq{bhyp}
 a(x)\geq  
a_0\qquad\mbox{a.e. }x\in\RR^N\quad \mbox{\ and }\quad 
{\rm ess}\lim_{|x|\to\infty}a(x)=+\infty\, .
\endeq

We assume throughout this paper that $p\in C_+(\RR^N)$ ($N\geq 2$) such that $p^+$ is finite.

Let
$f:\RR^N\times\RR\ri\RR$ be a  measurable function
such that, for some $C>0$,  $q\in\RR$ with  
$p^+<q+1\leq Np^-/(N-p^-)$ if $p^-<N$ and $p^+<q+1<+\infty$ if $p^-\geq N$, 
and $\mu > p^+$, we have 
\neweq{f1hyp}
|f(x,t)|\leq  C(|t|+|t|^q)\qquad\mbox{a.e. }(x,t)\in\RR^N\times\RR\,;
\endeq
\neweq{f2hyp}
\lim_{\ep\searrow 0}\ {\rm esssup}\ \bigg\{\bigg|\frac{f(x,t)}{t^{p^+-1}}\bigg|;\
(x,t)\in\RR^N\times(-\ep,\ep)\bigg\}=0\,;
\endeq
\neweq{f3hyp}
0\leq \mu F(x,t)\leq 
t\underline f (x,t)\quad\mbox{a.e. }(x,t)\in\RR^N\times[0,+\infty)\, .
\endeq

Our hypothesis $q\leq Np^-/(N-p^-)$ enables us to allow an almost critical behaviour on $f$.   
We also point out that we do not assume that the nonlinearity $f$
is continuous. An example of nonlinearity that fulfills assumptions \eq{f1hyp}--\eq{f3hyp}
is any positive discontinuous function $f(x,t)$ with subcritical growth 
that obeys like $t^a$ ($a>p^+$) in a neighborhood of the origin; for instance, take $N=3$,
$p^-=2$, $p^+=5$ (say, $p(x)=(7+3\,\sin |x|)/2$), and
$$f(x,t)=\left\{\begin{array}{lll}
&\di t^5+t^6&\di\qquad\mbox{if $(x,t)\in\RR^N\times [0,1)$}\\
&\di e^t&\di\qquad\mbox{if $(x,t)\in\RR^N\times [1,10]$}\\
&\di \rho (x)+t^5&\di\qquad\mbox{if $(x,t)\in\RR^N\times (10,\infty)$}\\
&\di -f(x,t)&\di\qquad\mbox{if $(x,t)\in\RR^N\times (-\infty,0)\,,$}\end{array}\right.$$
where $\rho(x)=+1$ if $|x|\in\QQ$ and $\rho (x)=0$, otherwise.

Let $E$ denote the set of all measurable functions $u:\RR^N\ri\RR$ such that 
$\left[a(x)\right]^{1/p(x)}u\in L^{p(x)}(\RR^N)$ and $|\nabla u|\in L^{p(x)}(\RR^N)$.
Then $E$ is a Banach space if it is endowed with the norm
 $$\|u\|_E:=\left| \left[a(x)\right]^{1/p(x)}\, u\right|_{p(x)}+
 |\nabla u|_{p(x)}\,.$$ We remark that $E$ is continuously embedded in $W^{1,p(x)}(\RR^N)$.
In the case $p(x)\equiv 2$ and if the potential $a(x)$ fulfills more general hypotheses
than \eq{bhyp}, then the embedding $E\subset L^{q+1}(\RR^N)$ is compact, whenever $2\leq q<
(N+2)/(N-2)$
(see, e.g., Bartsch, Liu and Weth \cite{blw} and 
Bartsch, Pankov and Wang \cite{bpw}). We do not know if this compact embedding 
still holds true in our ``variable exponent" framework and under assumption \eq{bhyp}.
 
Throughout this paper we denote by $\langle\cdot,\cdot\rangle$ the duality pairing between
$E^*$ and $E$.

Set
$F(x,t):=\int_0^tf(x,s)ds$ and
$$\Psi (u):=\int_{\RR^N} F(x,u(x))dx\,.$$ We observe that $\Psi$
is locally Lipschitz on $E$. This follows by \eq{f1hyp}, H\"{o}lder's
inequality and the continuous embedding $E\subset L^{q+1}(\RR^N)$. Indeed, for all $u,\,v\in E$,
$$
|\Psi(u)-\Psi(v)|\leq  C\, \|u-v\|_E\, ,
$$
where $C=C(\|u\|_E,\|v\|_E) > 0$ depends only on $\max\{\|u\|_E,\|v\|_E\}$.

In this paper we are concerned with the problem
\neweq{Phyp}\left\{\begin{array}{lll}
&\di -{\rm div}\, (|\nabla u|^{p(x)-2}\nabla u)+a(x)|u|^{p(x)-2}u\in[\underline f(x,u),\overline 
f(x,u)]\qquad &\di\mbox{in }\RR^N\\
&\di u\geq 0,\ u\not\equiv 0\qquad &\di\mbox{in }\RR^N\, .\end{array}\right.
\endeq

We notice that the semilinear anisotropic 
case corresponding to $p(x)\equiv 2$ has been analyzed in Gazzola and R\u adulescu \cite{6}.

We refer to Bertone--do \'O \cite{bertone} and Krist\'aly \cite{ksandor} for the 
study (by means of other methods) of certain classes of Schr\"odinger type equations which involve
discontinuous nonlinearities.

\begin{defin}\label{def1hyp}
We say that $u\in E$ is a solution of Problem \eq{Phyp} if $u\geq 0$, $\not\equiv 0$, and
$0\in\partial I(u)$,
where
$$ I(u):=\int_{\RR^N}\frac{1}{p(x)}\left(|\nabla u|^{p(x)}+a(x)|u|^{p(x)}\right)dx - \int_{\RR^N} 
F(x,u^+)dx\,,
\qquad\mbox{for all } u\in E\,.
$$
\end{defin}

The mapping $I:E\ri\RR$ is called the energy functional associated to Problem \eq{Phyp}. Our 
previous remarks show that $I$ is locally Lipschitz on the Banach space $E$.

The above definition may be reformulated, equivalently, in terms of hemivariational
inequalities. More precisely, $u\in E$ is a solution  of \eq{Phyp} if $u\geq 0$, $u\not\equiv 0$
in $\RR^N$, and
$$
\int_{\RR^N}\left(|\nabla u|^{p(x)-2}\nabla u\nabla 
v+a(x)|u|^{p(x)-2}uv\right)dx+\int_{\RR^N}(-F)^0(x,u;v)dx\geq  0,\qquad\mbox{for all } v\in E\,.
$$

Our main result is the following

\begin{teo}\label{t1hyp}
Assume that hypotheses \eq{bhyp}--\eq{f3hyp} are fulfilled. Then Problem \eq{Phyp}
has at least one solution.
\end{teo}

\section{Proof of Theorem \ref{t1hyp}}
We first claim that
there exist positive constants $C_1$ and $C_2$  such that
\neweq{4}
f(x,t)\geq   C_1t^{\mu-1}-C_2\qquad\mbox{a.e. }(x,t)\in\RR^N
\times[0,+\infty)\, .
\endeq
Indeed, by the definition of $\underline f$ we deduce that
\neweq{5}
\underline f(x,t)\leq  f(x,t)\qquad\mbox{a.e. }(x,t)\in\RR^N\times[0,+\infty)\, .
\endeq
Set
$\underline F(x,t):=\int_0^t\underline f(x,s)ds$.
Thus, by our assumption \eq{f3hyp},
\neweq{6}
0\leq \mu\underline F(x,t)\leq  t\underline f(x,t)\qquad\mbox{a.e. }(x,t)\in
\RR^N\times[0,+\infty)\, .
\endeq
Next, by \eq{6}, there exist  positive constants $R$ and $K_1 $ such that
\neweq{7}
\underline F(x,t)\geq   K_1t^\mu\qquad\mbox{a.e. }(x,t)\in\RR^N
\times[R,+\infty )\, .
\endeq
Our claim \eq{4} follows now directly by relations \eq{5}, \eq{6} and \eq{7}.

Next, we observe that
$$\partial I(u)=-{\rm div}\,(|\nabla u|^{p(x)-2}\nabla u) 
+a(x)|u|^{p(x)-2}u-\partial\Psi(u^+)\qquad\mbox{in }E^*\,.$$
So, by  \cite[Theorem 2.2]{2} and  \cite[Theorem 3]{mr}, we have
$$\partial\Psi (u)\subset [\underline f(x,u(x)),\overline f(x,u(x))]\qquad
\mbox{a.e. }x\in\RR^N\, ,$$
in the sense that if $w\in\partial\Psi(u)$ then
\neweq{w}
\underline f(x,u(x))\leq  w(x)\leq \overline f(x,u(x))\qquad
\mbox{a.e. }x\in\RR^N\, .
\endeq
This means that if $u_0$ is a critical point of $I$, then there exists
$w\in\partial\Psi(u_0)$ such that
$$-{\rm div}\,(|\nabla u_0|^{p(x)-2}\nabla u_0) +a(x)|u_0|^{p(x)-2}u_0=w\qquad\mbox{in} \ E^*\, 
.$$
This argument shows that, for proving Theorem \ref{t1hyp}, it is enough to show that the energy 
functional $I$ has at least a nontrivial critical point $u_0\in E$, $u_0\geq 0$.
We prove the existence of a solution of Problem \eq{Phyp} by arguing that the hypotheses of 
Chang's version of the
Mountain Pass Lemma for locally Lipschitz functionals (see Chang \cite{2}) are fulfilled. 
  More precisely, we check the following geometric  assumptions:
\neweq{h1}
I(0)=0\ \mbox{and there exists}\  v\in E\ \mbox{such that}\ I(v)\leq  0\,;
\endeq
\neweq{h2}
\mbox{there exist }\beta,\rho > 0\quad\mbox{such that}\quad I\geq  \beta\quad\mbox{on}\quad
\{u\in E;\ \|u\|_E=\rho\}\, .
\endeq

\smallskip
{\sc Verification of \eq{h1}.} 
Fix
$w\in\C(\RR^N)\setminus\{0\}$ such that $w\geq  0$ in $\RR^N$. In particular, we have
$$\int_{\RR^N}\left(|\nabla w|^{p(x)}+a(x)w^{p(x)}\right)dx < +\infty\, .$$
So, by \eq{4} and choosing $t>1$,
$$\begin{array} {ll}
\di I(tw)&\di =\int_{\RR^N}\frac{t^{p(x)}}{p(x)} \left(| \nabla w|^{p(x)}+a(x)w^{p(x)}\right)dx- 
\Psi (tw) \\ & \di\leq \frac{t^{p^+}}{p_-}\int_{\RR^N}\left(|\nabla 
w|^{p(x)}+a(x)w^{p(x)}\right)dx+C_2t\, \int_{\RR^N} wdx-C_1't^\mu\
\int_{\RR^N} w^\mu dx \,.
\end{array}
$$
Since, by hypothesis, $1< p^+<\mu$, we deduce that $I(tw)<0$ for $t > 1$ large enough.

\smallskip
{\sc Verification of \eq{h2}.} Our hypotheses \eq{f1hyp} and \eq{f2hyp} imply
that, for any $\ep > 0$, there exists some $C_\ep >0$ such that
\neweq{growth}
|f(x,t)|\leq \ep|t|+C_\ep|t|^q\qquad\mbox{ a.e. }(x,t)\in\RR^N\times\RR\, .
\endeq
By \eq{growth} and Sobolev embeddings in variable exponent spaces we have, for any $u\in E$,
$$
\Psi(u)\leq \ep \int_{\RR^N}\frac{1}{p(x)} |u|^{p(x)}dx+\frac{A_\ep}{q+1}\ 
\int_{\RR^N}|u|^{q+1}dx\leq 
\ep  \int_{\RR^N}\frac{1}{p(x)} |u|^{p(x)}dx+C_4\, \|u\|_{L^{q+1}(\RR^N)}^{q+1},
$$
where $\eps$ is arbitrary and $C_4=C_4(\eps)$. Thus, by our hypotheses,
$$\begin{array}{ll}
\di I(u)&\di =\int_{\RR^N}\frac{1}{p(x)}\left(|\nabla 
u|^{p(x)}+a(x)|u|^{p(x)}\right)dx-\Psi(u^+)\\
&\di\geq  \frac{1}{p^+}
\int_{\RR^N}\left[|\nabla u|^{p(x)}+(a_0-\ep)|u|^{p(x)}\right]dx
 -C_4\, \|u\|_{L^{q+1}(\RR^N)}^{q+1}\geq  \beta > 0\, ,\end{array}
$$
for $\|u\|_E=\rho$, with $\rho$, $\ep$ and $\beta$ are small enough positive
constants.\qed

Denote
$${\mathcal P}:=\{\gamma\in C([0,1],E);\ \gamma (0)=0,\ \gamma(1)\not=0
\mbox{ and }I(\gamma(1))\leq  0\}$$
and
$$c:=\inf_{\gamma\in{\mathcal P}}\,\max_{t\in [0,1]}\, I(\gamma (t))\, .$$
Set
$$\lambda_I(u):=\min_{\zeta\in\partial I(u)}\|\zeta\|_{E^*}\, .$$
We are now in position to apply Chang's version of the Mountain Pass Lemma for locally Lipschitz 
functionals (see Chang \cite{2}). So, there exists a sequence
$\{u_n\}\subset E$ such that
\neweq{palais}
I(u_n)\ri c\qquad\mbox{and}\qquad\lambda_I(u_n)\ri 0\,.
\endeq
Moreover, since $I(|u|)\leq  I(u)$ for all $u\in E$, we can assume without loss of generality that
$u_n\geq 0$ for every $n\geq 1$. So, for all positive integer $n$, there exists 
$\{w_n\}\in\partial\Psi(u_n)\subset E^*$ such that, for any $v\in E$,
\neweq{psi1}
\int_{\RR^N}\left(|\nabla u_n|^{p(x)-2}\nabla u_n\nabla v +a(x)u_n^{p(x)-1}v\right)dx-\langle 
w_n,v\rangle\ri 0\qquad\mbox{as $n\ri\infty$}\, .
\endeq
Note that for all $u\in E$, $u\geq 0$, the definition of $\Psi$ and our hypotheses yield
$$\Psi(u)\leq \frac{1}{\mu}\ \int_{\RR^N} u(x)\underline f(x,u(x))dx\, .$$
Therefore, by \eq{w}, for every $u\in E$, $u\geq 0$, and for any $w\in\partial\Psi(u)$,
$$
\Psi(u)\leq \frac{1}{\mu}\ \int_{\RR^N} u(x)w(x)dx\, .
$$
Hence
$$\begin{array}{ll}
I(u_n) &\geq \di\frac{\mu -p^+}{\mu p^+}\, \int_{\RR^N}\left(|\nabla 
u_n|^{p(x)}+a(x)u_n^{p(x)}\right)dx\\
 &+\di\frac{1}{\mu} \int_{\RR^N}\left(|\nabla u_n|^{p(x)}+a(x)u_n^{p(x)}-w_nu_n\right)dx
+\frac{1}{\mu}\int_{\RR^N}w_nu_ndx-\Psi(u_n)\\
&\geq   \di\frac{\mu -p^+}{\mu p^+}\, \int_{\RR^N}\left(|\nabla 
u_n|^{p(x)}+a(x)u_n^{p(x)}\right)dx
+\di\frac{1}{\mu} \int_{\RR^N}\left(|\nabla u_n|^{p(x)}+a(x)u_n^{p(x)}-w_nu_n\right)dx\\
&=\di\frac{\mu -p^+}{\mu p^+}\, \int_{\RR^N}\left(|\nabla u_n|^{p(x)}+a(x)u_n^{p(x)}\right)dx
+\di\frac{1}{\mu}\langle -\Delta_{p(x)}u_n+au_n-w_n,u_n\rangle\\
&=\di\frac{\mu -p^+}{\mu p^+}\, \int_{\RR^N}\left(|\nabla u_n|^{p(x)}+a(x)u_n^{p(x)}\right)dx+ 
o(1)\| u_n\|_{E}\, .
\end{array}$$
This relation and \eq{palais} show that the Palais-Smale sequence
$\{u_n\}$ is bounded in $E$. It follows that $\{u_n\}$ converges weakly (up to a subsequence)
in $E$ and strongly in $L^{p(x)}_{\rm loc}(\RR^N)$ to some $u_0\geq 0$. Taking into
account that $w_n\in\partial\Psi(u_n)$ for all $m$, that $u_n\weak u_0$ in $E$
and that there exists $w_0\in E^*$ such that $w_n\weak w_0$ in $E^*$ (up to a
subsequence), we infer that $w_0\in\partial\Psi(u_0)$. This follows from the
fact that the map $u\longmapsto F(x,u)$ is compact from $E$ into $L^1$. Moreover,
if we take $\phi\in C^\infty_{\rm c}(\RR^N)$ and let $\omega:=\supp\,\phi$,
then by \eq{psi1} we get
$$\int_\omega\left(|\nabla u_0|^{p(x)-2}\nabla 
u_0\nabla\varphi+a(x)u_0^{p(x)-1}\varphi-w_0\varphi\right)dx=0\,. $$
So, by relation (4) p.104 in Chang \cite{2} and by the
definition of $(-F)^0$, we deduce that
$$\int_\omega\langle(|\nabla u_0|^{p(x)-2}\nabla 
u_0\nabla\varphi+a(x)u_0^{p(x)-1}\varphi\rangle)dx+\int_\omega
(-F)^0(x,u_0;\varphi)dx\geq  0\, .$$
By density, this hemivariational inequality holds for all $\phi\in E$ and
this means that $u_0$ solves Problem \eq{Phyp}.

It remains to prove that $u_0\not\equiv 0$. If $w_n$ is as in \eq{psi1}, then
by \eq{w} (recall that $u_n\geq 0$) and \eq{palais} (for large $m$) we deduce that
\neweq{c2}\begin{array}{ll}
\di\frac c2&\di\leq  I(u_n)-\frac{1}{p^-}\, \langle-\Delta_{p(x)} u_n+au_n-w_n,u_n\rangle\\
&\di=\frac{1}{p^-}\, \langle w_n,u_n\rangle-\int_{\RR^N} F(x,u_n)dx
\leq  \frac{1}{p^-}\,\int_{\RR^N} u_n\overline f(x,u_n)dx\, .\end{array}
\endeq
Now, taking into account its definition, one deduces that $\overline f$
verifies \eq{growth}, too. So, by \eq{c2}, we obtain
$$
0<\frac{c}{2}\leq \frac{1}{p^-}\int_{\RR^N}(\ep u_n^2+A_\ep u_n^{q+1})dx
=\frac{\ep}{p^-}\, \|u_n\|^2_{L^2(\RR^N)}+\frac{A_\ep}{p^-}\, \|u_n\|^{q+1}_{L^{q+1}(\RR^N)}\,.
$$
In particular, this shows that $\{u_n\}$ does not converge strongly to 0 in $L^{q+1}(\RR^N)$. 
It remains to argue that $u_0\not\equiv 0$.
Since both $\|u_n\|_{L^{p^-}(\RR^N)}$ and $\|\nabla u_n\|_{L^{p^-}(\RR^N)}$
are bounded, it follows by Lemma I.1  in Lions \cite{pll} that the sequence
$\{u_n\}$ ``does not vanish'' in $L^{p^-}(\RR^N)$. Thus, there exists a sequence
$\{z_n\}\subset\RR^N$ and $C > 0$ such that, for some $R>0$,
\begin{equation}\label{nonvan}
\int_{z_n+B_R}u_n^{p^-}dx\geq   C\,.
\end{equation}
 We claim that the sequence $\{z_n\}$ is bounded in $\RR^N$. Indeed,  if not, up to
a subsequence, it follows by \eq{bhyp} that
$$\int_{\RR^N} a(x)u_n^{p^-}dx\ri +\infty\qquad\mbox{as $n\ri\infty$}\,,$$
which contradicts our assumption $I(u_n)=c+o(1)$. Therefore, by \eq{nonvan}, there
exists an open bounded set $D\subset\RR^N$ such that
$$
\int_D u_n^{p^-}dx\geq   C > 0\, .
$$
In particular, this relation implies that $u_0\not\equiv 0$ and our proof is concluded.
\qed

\medskip
{\bf Acknowledgments}. The author is grateful to the anonymous referee for the
careful reading of the manuscript and for several useful remarks. This work
is a part of the author's Ph.D. thesis at the Babe\c s--Bolyai University in Cluj.
I am very pleased to acknowledge my adviser, Professor Radu Precup, for his constant support
and high level guidance  during the preparation of this thesis.

\end{document}